# A hybrid extended finite element/level set method for modeling equilibrium shapes of nano-inhomogeneities

**Xujun Zhao[1], Stéphane P.A. Bordas[3] and Jianmin Qu[1,2]**

[1]Department of Mechanical Engineering, Northwestern University, Evanston, IL USA 60208
[2]Department of Civil and Environmental Engineering, Northwestern University, Evanston, IL USA 60208
[3]Cardiff School of Engineering, Institute of Mechanics and Advanced Materials Theoretical, Applied and Computational Mechanics, Cardiff University, The Parade, Cardiff CF24 3AA, Wales, U.K.

E-mail: j-qu@northwestern.edu, stephane.bordas@alum.northwestern.edu.

**Abstract**
Interfacial energy plays an important role in equilibrium morphologies of nanosized microstructures of solid materials due to the high interface-to-volume ratio, and can no longer be neglected as it does in conventional mechanics analysis. The present work develops an effective numerical approach by means of a hybrid smoothed extended finite element/level set method to model nanoscale inhomogeneities with interfacial energy effect, in which the finite element mesh can be completely independent of the interface geometry. The Gurtin-Murdoch surface elasticity model is used to account for the interface stress effect and the Wachspress interpolants are used for the first time to construct the shape functions in the smoothed extended finite element method. Selected numerical results are presented to study the accuracy and efficiency of the proposed method as well as the equilibrium shapes of misfit particles in elastic solids. The presented results compare very well with those obtained from theoretical solutions and experimental observations, and the computational efficiency of the method is shown to be superior to that of its most advanced competitor.

## 1. Introduction

Recent advances in nanotechnology exacerbate the need for computational tools that are capable of capturing the effects of interfaces, which play an important role in nanostructured materials due to a characteristically high interface-to-volume ratio. Although atomic level computational tools such as molecular dynamic (MD) and first principle calculations are able to simulate interface effects, these methods are computationally intensive, thus their applications are usually limited to nano-scale samples and nano-second time durations. Many engineering problems, however, occur at much larger spatial and temporal scales. For example, simulating the formation and morphological evolution of precipitates in superalloys involves length scales ranging from several nanometers to tens of micrometers, and the physical processes last up to hours in time. In these cases, it is necessary to use a continuum level model that can capture the interfacial effects of particles at different length and time scales.

The bonding environment of atoms near surfaces/interfaces differs from that in the bulk, so that the energy associated with surfaces/interfaces is also different from the bulk energy. This effect will be prominent for nanostructured materials when the surface-to-volume ratio is high. Gurtin and Murdoch [16, 17] developed a generic continuum model incorporating the surface/interface effects where the surface/interface is modeled as a zeor-thickness layer with its own physical properties. This surface/interface elasticity model has been extensively used to model nanostructures [8, 24, 33, 43, 7], and it was shown that the size-dependent behaviors obtained from the surface elasticity model agree very well with the atomistic simulations [28, 24].

However, these theoretical results are limited to systems with either simple geometries or isotropic material properties. To characterize the size dependent behavior of nanosized structures with more realistic geometries and anisotropic properties, finite element methods (FEM) incorporating interface effects were developed by introducing surface/interface elements [14, 40, 34]. Within the standard FEM framework, it is required that both the bulk mesh and the interface mesh conform to the material interface in order to describe the strain discontinuities across the interface, so that remeshing is generally required when the interface moves. This usually leads to difficulties for inhomogeneities with very complex and evolving geometries, especially when the inhomogeneities are numerous and evolving with time.

The extended finite element method (XFEM), which was originally developed to model cracks in fracture mechanics [26], is a robust and powerful computational tool capable of modeling arbitrary discontinuities without requiring the meshes to be conformed with the discontinuity interfaces. The XFEM has been successfully used to solve many physical problems with arbitrary interfaces, such as holes and inclusions [37], phase solidification [5], multiphase flows [4], biofilm growth [9, 10], etc. A comprehensive review of XFEM can be found in [13] and an open source C++ XFEM code is also available for use [1]. Yvonnet et al. [42] and Farsad et al. [12] have recently presented an XFEM approach to study the surface effects on the mechanical behavior of nano-scale materials without considering the moving interface and the dynamic level set.

In this paper, we developed a hybrid smoothed extended finite element/level set method incorporating the interfacial energy effect, in which a level set function is used to implicitly capture the motion of the interface and the finite element mesh is independent of the interface geometry, so that no remeshing is required even when the interfaces move, merge, cusp or disappear. The strain smoothing technique [21] is applied to modify the standard XFEM so that the domain integrals in each element can be transformed into contour integrals along the boundaries of smoothing cells. Such a transformation alleviates the needs of numerically evaluating the Jacobian matrix and derivatives of the shape functions, which reduces the computational cost without loss of accuracy. The Wachspress interpolants [41] are used to construct the shape functions, which greatly facilitate the numerical quadrature along the boundaries of smoothing cells. We present some numerical examples to illustrate the computational accuracy and convergence of the Smoothed XFEM, followed by an example to study the influence of interface effects on the equilibrium shape of particles with different size.

## 2. Formulation

Consider a finite elastic domain $V$ containing multiple inhomogeneities with arbitrary shapes and misfit strain $\boldsymbol{\varepsilon}^*$ (eigenstrain). $\Omega$ is the union of all the inhomogeneities and $\Gamma = \partial\Omega$ denotes the interfaces between the inhomogeneities and the matrix. Further, we assume that $V$ is subjected to a prescribed body force **b**, a surface traction **t** on $S_t$ and a displacement $\mathbf{u}_0$ on $S_u$, where $S = S_u \bigcup S_t$.

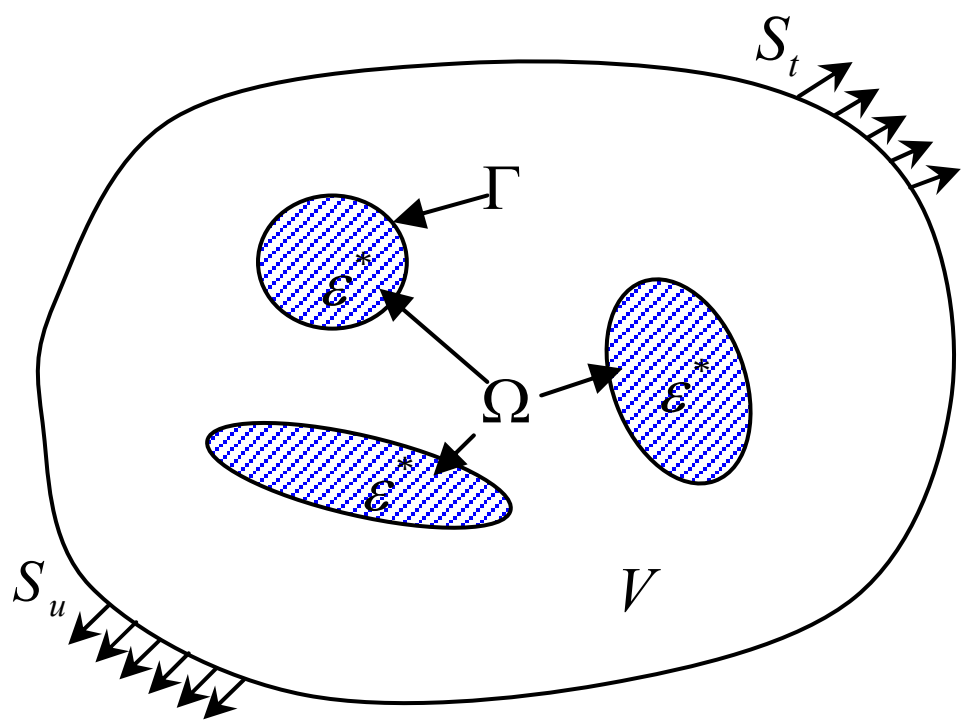

Figure 1. Inhomogeneities in a finite matrix

In terms of the displacement vector **u** and the total strain tensor $\boldsymbol{\varepsilon}$, the total potential energy $\Pi$ of the composite consists of three parts: the elastic energy $U^B$ in the bulk, the interfacial energy $U^s$ on all the inhomogeneity-matrix interfaces, and the work done by external forces *W*, i.e.,

$$\Pi = U^B + U^S + W \tag{1}$$

where

$$U^B = \frac{1}{2}\int_{V\backslash\Omega} \boldsymbol{\varepsilon} : \mathbf{L}^M : \boldsymbol{\varepsilon} \; d\Omega + \frac{1}{2}\int_{\Omega} (\boldsymbol{\varepsilon} - \boldsymbol{\varepsilon}^*) : \mathbf{L}^I : (\boldsymbol{\varepsilon} - \boldsymbol{\varepsilon}^*) \; d\Omega \tag{2}$$

$$U^S = \int_{\Gamma} \gamma \; dS \tag{3}$$

$$W = -\int_V \mathbf{u} \cdot \mathbf{b} \; d\Omega - \int_{S_t} \mathbf{u} \cdot \mathbf{t} \; dS \tag{4}$$

where $\gamma$ is the interfacial excess energy density [7] given by

$$\gamma = \gamma_0 + \boldsymbol{\tau}^S : \boldsymbol{\varepsilon}^S + \frac{1}{2}\boldsymbol{\varepsilon}^S : \mathbf{L}^S : \boldsymbol{\varepsilon}^S \tag{5}$$

In the above, field quantities with a superscript *I*, *M* and *S* are associated with the inhomogeneity, the matrix and the interfaces, respectively. For example, $\mathbf{L}^I$ and $\mathbf{L}^M$ represent the elastic stiffness tensors of inhomogeneities and matrix. $\boldsymbol{\varepsilon}^S$ is the interface strain tensor, and $\mathbf{L}^S$ is the interface elastic stiffness tensor. $\boldsymbol{\tau}^S$ is the interface residual stress and $\gamma_0$ is the interfacial excess energy density that corresponds to the interfacial energy state when $\boldsymbol{\varepsilon}^S = \mathbf{0}$. It should be mentioned that the constant $\gamma_0$ has no influence on the elastic field for a specific configuration but contributes to the total system energy.

The corresponding interfacial stress can be obtained from the Shuttleworth equation [35],

$$\boldsymbol{\sigma}^S = \frac{\partial \gamma}{\partial \boldsymbol{\varepsilon}^S} = \boldsymbol{\tau}^S + \mathbf{L}^S : \boldsymbol{\varepsilon}^S \tag{6}$$

where $\mathbf{I}$ is the second order identity tensor. According to the coherent interface assumption, the interfacial strain can be obtained by taking the tangential gradient of the displacement vector [16]

$$\boldsymbol{\varepsilon}^S = \frac{1}{2}\left[\nabla_S \otimes \mathbf{u} + \left(\nabla_S \otimes \mathbf{u}\right)^T\right] \tag{7}$$

where $\nabla_S = \mathbf{P} \cdot \nabla$ is the interfacial gradient, $\mathbf{P} = \mathbf{I} - \mathbf{n} \otimes \mathbf{n}$ is the tangential projection operator, and $\mathbf{n}$ is the outward unit normal to the interface.

Under equilibrium, the system minimizes its potential energy and the stationary condition requires the variation of the functional to vanish, i.e. $\delta\Pi = 0$, which gives

$$\delta\Pi = \delta U^b + \delta U^S + \delta W = 0 \tag{8}$$

where

$$\begin{aligned}
\delta U^b &= \int_{V-\Omega} \boldsymbol{\varepsilon} : \mathbf{L}^M : \delta\boldsymbol{\varepsilon}\, d\Omega + \int_\Omega (\boldsymbol{\varepsilon} - \boldsymbol{\varepsilon}^*) : \mathbf{L}^I : \delta\boldsymbol{\varepsilon}\, d\Omega \\
\delta U^S &= \int_\Gamma [\boldsymbol{\tau}^S : \delta\boldsymbol{\varepsilon}^S + \boldsymbol{\varepsilon}^S : \mathbf{L}^S : \delta\boldsymbol{\varepsilon}^S]\, dS \\
\delta W &= -\int_\Omega \mathbf{b} \cdot \delta\mathbf{u}\, d\Omega - \int_{S_t} \mathbf{t} \cdot \delta\mathbf{u}\, dS
\end{aligned} \tag{9}$$

Equation (8) will be used later to develop the weak form of the finite element equations.

## 3. Level set description of the interfaces

The level set method was originally devised for tracking a moving interface [30], and became a key ingredient of XFEM to implicitly describe complicated geometrical interfaces of microstructures without tracking them explicitly, such as cracks [38], holes and inclusions [37, 27], dislocations [15] and biofilms [9, 10], as well as fluid-structure interfaces [20].

In this paper, the interfaces between the inhomogeneities and the matrix can be defined as the zero level set of a function $\phi(\mathbf{x})$, such that $\phi < 0$ in the particles, $\phi > 0$ in the matrix and $\phi = 0$ on the interface, which is as shown in Figure 2. Often, the signed distance is used as

$$\phi(\mathbf{x}) = \|\mathbf{x} - \mathbf{x}_{\min}\| sign(\mathbf{n}_{\min} \cdot (\mathbf{x} - \mathbf{x}_{\min})) \tag{10}$$

where $\mathbf{x}_{\min}$ is the orthogonal projection of $\mathbf{x}$ on the interface $\Gamma$, $\mathbf{n}_{\min}$ is the outward unit normal at $\mathbf{x}_{\min}$, and $sign(\bullet)$ denotes the sign function.

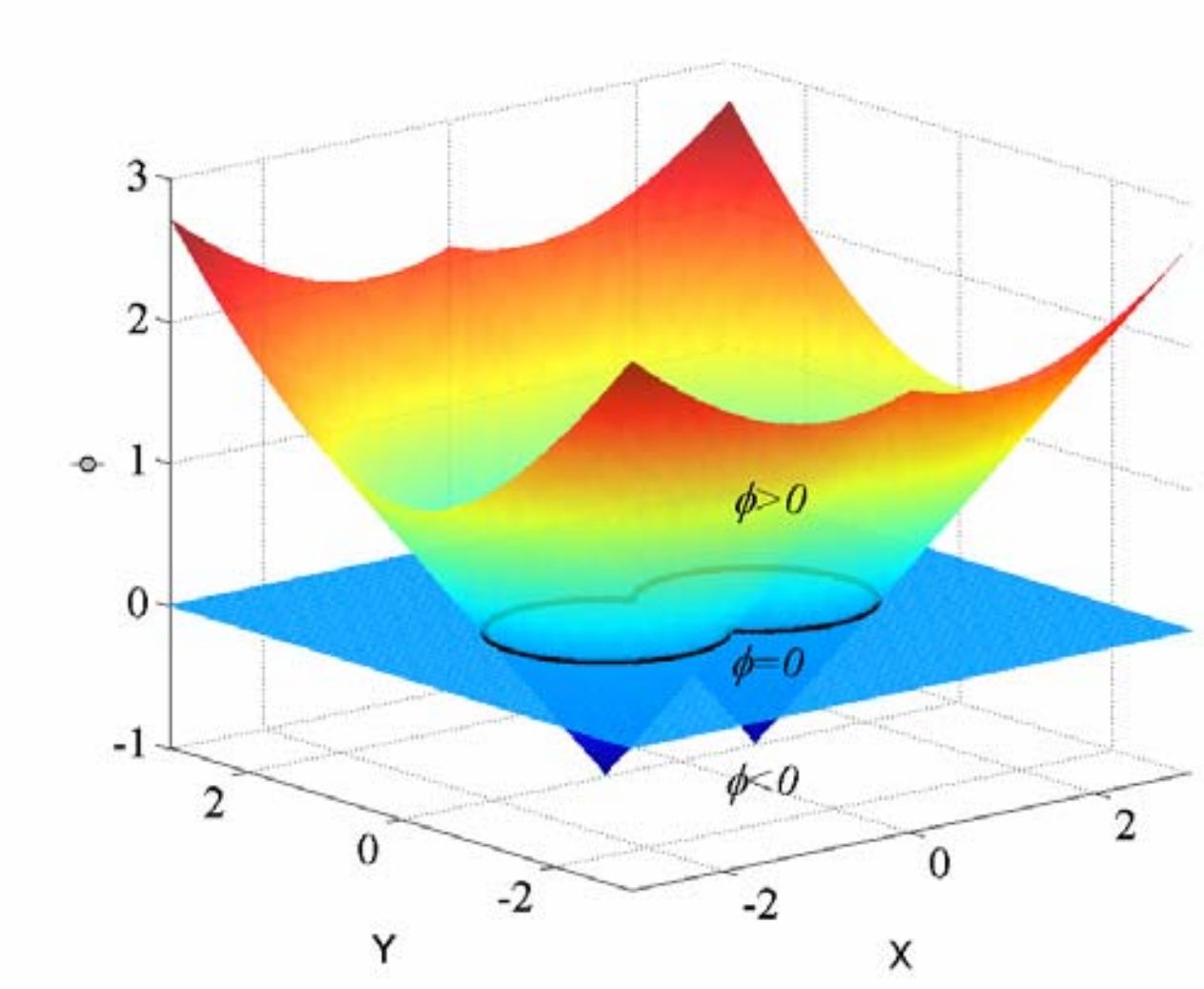

Figure 2. Level set function

It should be pointed out that the function $\phi(\mathbf{x})$ is usually not known explicitly, or analytically, except for simple shapes (circles, ellipses, etc.), but can be approximated by a set of local values, $\phi_I$, at the finite element nodes:

$$\phi(\mathbf{x}) = \sum_{I=1}^{M} N_I(\mathbf{x}) \phi_I \tag{11}$$

where $N_I(\mathbf{x})$ is the finite element shape function of node *I*, and *M* is the total number of nodes in an element.

An initial value partial differential equation then can be obtained for the evolution of $\phi$ by taking a material derivative on both sides of (10) [30]

$$\frac{\partial \phi}{\partial t} + v_n^{ext} \left|\nabla \phi\right| = 0 \tag{12}$$

which is the well-known level set equation. The position of the interface at time $t$ is defined by the zero level set $\phi(\mathbf{x},t)=0$. $v_n^{ext} = \mathbf{v}\cdot\mathbf{n}$ is the normal velocity on the interface

A main advantage of the level set description is its ability to describe an arbitrary number of inhomogeneities with a single level set function. The geometric quantities such as normal vector $\mathbf{n}$ and curvature $\kappa$ can be easily expressed in terms of level set function, e.g.

$$\mathbf{n} = \frac{\nabla\phi}{\left|\nabla\phi\right|}, \qquad \kappa = \nabla\cdot\frac{\nabla\phi}{\left|\nabla\phi\right|} \tag{13}$$

Additionally, when the interfaces evolve in time, merge, create cusps, coalesce, the hybrid smoothed XFEM/level set method can deal with the topological transformation naturally without remeshing or other special treatments.

## 4. Smoothed extended finite element method

In this section, we present the basic idea of the strain smoothing technique and formulate the Smoothed XFEM (SmXFEM) to account for the effects of interfaces. This will be based on the Wachspress shape functions[41]. The implementation of boundary integration and smoothing cell subdivision will also be discussed.

### 4.1 *Strain smoothing in finite element*

The strain smoothing technique was first proposed in the context of mesh-free methods to stabilize nodal integration [3]. Liu et al. [21] introduced this method into the conventional finite element formulation. Without introducing additional degrees of freedom, an element is further subdivided into several smoothing cells, and a smoothing operation is performed in each cell within an element by a weighted average of the standard FEM strain field in the Voigt notation , $\boldsymbol{\varepsilon}^h \equiv [\varepsilon_{11}^h, \varepsilon_{22}^h, \gamma_{12}^h]^T$. We note that the Voigt notation for the strain tensors will be used in the rest of this paper, except Section 5 where the tensorial notation is used.

For example, the smoothed strain value at a point $\mathbf{x}_k$ can be expressed as

$$\bar{\boldsymbol{\varepsilon}}(\mathbf{x}_k) = \int_{\Omega_k^s} \boldsymbol{\varepsilon}^h(\mathbf{x}) W(\mathbf{x}_k - \mathbf{x}) d\Omega \tag{14}$$

where $\Omega_k^s$ is the smoothing cell defined in the local vicinity of $\mathbf{x}_k$ and $W(\mathbf{x}_k-\mathbf{x})$ is the smoothing or weight function associated with $\mathbf{x}_k$. The following piecewise constant weight function is generally used in the smoothed FEM formulation,

$$W(\mathbf{x}_k - \mathbf{x}) = \begin{cases} 1/A_k^s & \mathbf{x}\in\Omega_k^s \\ 0 & \mathbf{x}\notin\Omega_k^s \end{cases} \tag{15}$$

where $A_k^s$ is the area of the smoothing cell $\Omega_k^s$. Substituting equation (15) into equation (14), one can obtain the smoothed strain by application of the Green's theorem

$$\bar{\boldsymbol{\varepsilon}}(\mathbf{x}_k) = \frac{1}{A_k^s}\int_{\Gamma_k^s} \mathbf{L}_k(\mathbf{x})\cdot\mathbf{u}(\mathbf{x}) d\Gamma \tag{16}$$

where $\mathbf{L}_n(\mathbf{x})$ is a matrix function of the outward normal vector on the boundary $\Gamma_k^s$ of smoothing cell $\Omega_k^s$ and has the form

$$\mathbf{L}_k(\mathbf{x}) = \begin{bmatrix} n_x & 0 \\ 0 & n_y \\ n_y & n_x \end{bmatrix} \qquad (17)$$

In such a way, the domain integration over $\Omega_k^s$ becomes a boundary integration along the edges of the smoothing cells. An isoparametric mapping is not necessary and the evaluation of the Jacobian in each element can be avoided. Also, we can see that the displacement gradient does not appear in equation (16). As a result, the derivatives of the shape functions are not required in the smoothed finite element method.

4.2 *Extended finite element discretization*

For a coherent interface, the displacement is assumed to be continuous across the interface. The strain field, however, may not be continuous, which is commonly called *weak discontinuity*. In the standard finite element formulation, the element edges are required to coincide with the geometrical interfaces to guarantee the accuracy and optimal convergence. In contrast, XFEM allows the mesh to be almost independent of interfaces by introducing additional enrichment, so that an inhomogeneity can be modeled by simply modifying the level set function without remeshing when its geometry changes. The displacement can be approximated in the following way [37]

$$\mathbf{u}^h(\mathbf{x}) = \mathbf{N}(\mathbf{x})\mathbf{d} + \mathbf{N}_a(\mathbf{x})\mathbf{a} \qquad (18)$$

where $\mathbf{d}$ and $\mathbf{a}$ correspond to the nodal displacement vector and additional degrees of freedom due to enrichment, respectively, and $\mathbf{N}(\mathbf{x})$ and $\mathbf{N}_a(\mathbf{x})$ are the standard and enrichment shape function matrices

$$\mathbf{N}(\mathbf{x}) = \begin{bmatrix} \mathbf{N}_1 & \mathbf{N}_2 & \cdots & \mathbf{N}_M \end{bmatrix}, \quad \mathbf{N}_a(\mathbf{x}) = F(\mathbf{x})\mathbf{N} \qquad (19)$$

where

$$\mathbf{N}_I = \begin{bmatrix} N_I & 0 \\ 0 & N_I \end{bmatrix} \qquad (20)$$

The enrichment function $F(\mathbf{x})$ is given as suggested by [25]

$$F(\mathbf{x}) = \sum_{I=1}^{M} N_I(\mathbf{x})|\phi_I| - \left|\sum_{I=1}^{M} N_I(\mathbf{x})\phi_I\right| \qquad (21)$$

Note that the enrichment function is zero in the elements that do not contain any part of the interface, so there are no spurious terms in the blending elements as shown in Figure 3 [25].

Substituting equation (18) into (16), the corresponding smoothing strain can be readily obtained in the following matrix form

$$\bar{\boldsymbol{\varepsilon}} = \bar{\mathbf{B}}\mathbf{d} + \bar{\mathbf{B}}_a\mathbf{a} \qquad (22)$$

where

$$\bar{\mathbf{B}} = \frac{1}{A_k^s}\int_{\Gamma_k^s} \mathbf{L}_n(\mathbf{x})\mathbf{N}(\mathbf{x})d\Gamma, \qquad \bar{\mathbf{B}}_a = \frac{1}{A_k^s}\int_{\Gamma_k^s} \mathbf{L}_n(\mathbf{x})\mathbf{N}_a(\mathbf{x})d\Gamma \qquad (23)$$

On substituting displacement approximation (18) and smoothing strain (22) into the weak form, equations (8) and (9), a discrete linear system of equations can be easily obtained.

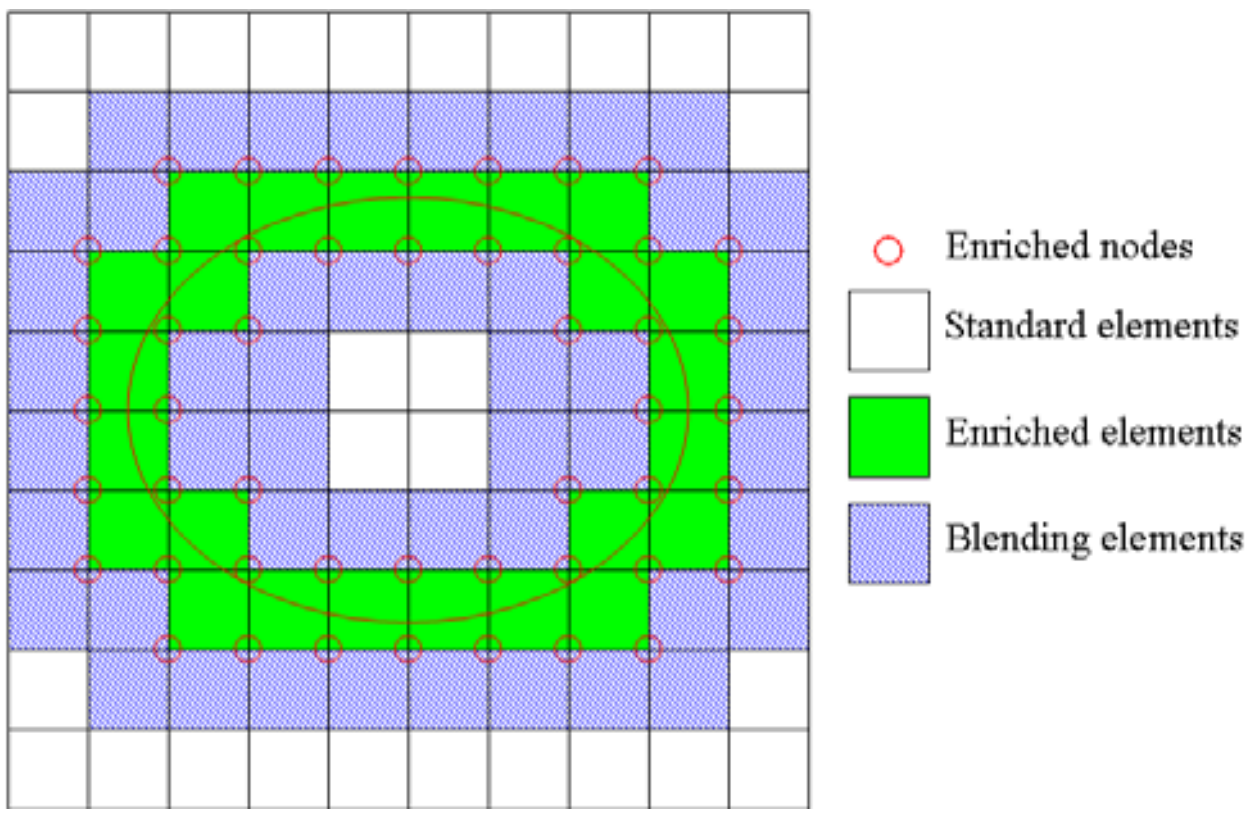

Figure 3. Enriched nodes (marked by red circles) which have additional degrees of freedom due to enrichment; enriched elements (marked by green squares) cut by the interface whose nodes are all enriched; blending elements (marked by shading squares) in which only part of the nodes are enriched.

*4.3 Wachspress shape functions and smoothing cell partition*

In the original paper on the smoothed FEM [21], non-mapped Lagrange shape functions are used to calculate the shape function values within a smoothed finite element, but the values at the vertices and mid-points of cell edges are evaluated using linear interpolations. This numerical scheme is valid only when the shape function is linear along the cell edges [2]. Unfortunately, this is not always the case. When the shape function is not linear along the edges in higher order elements or enriched elements in XFEM whose enrichment functions are not linear, the average shape function approximation and numerical integration scheme are not accurate. Bordas and Natarajan [2] suggested that the Wachspress interpolation can be used to retain the desirable features in the conventional smoothed FEM.

Following their ideas, we for the first time introduce the Wachspress shape functions into the SmXFEM framework to facilitate the accurate numerical integrals in the enriched finite elements. In this way, the shape function values at every edge of the smoothing cells can be directly obtained without any approximation. Moreover, since the Wachspress interpolants may be constructed on elements with arbitrary (including curved) edges, this makes our formulation general and attractive for curved interfaces. It is further shown that the numerical accuracy depends on the smoothing cell division used in the SmXFEM.

Although the Wachspress interpolants can be built for arbitrary *n*-sided polygonal elements with arbitrary edges, only quadrilateral elements shown in Figure 4 are formulated in the present work. Let $l_i(x, y) = 0$ be the function defining the line associated with the *i*-th edge of element $\Omega_e$ which can be uniquely written with coefficients $a_i$, $b_i$ and $c_i$ in the form,

$$l_i(x, y) = c_i - a_i x - b_i y \quad \text{for} \quad (x, y) \in \Omega_e \Rightarrow l_i(x, y) > 0 \quad . \tag{24}$$

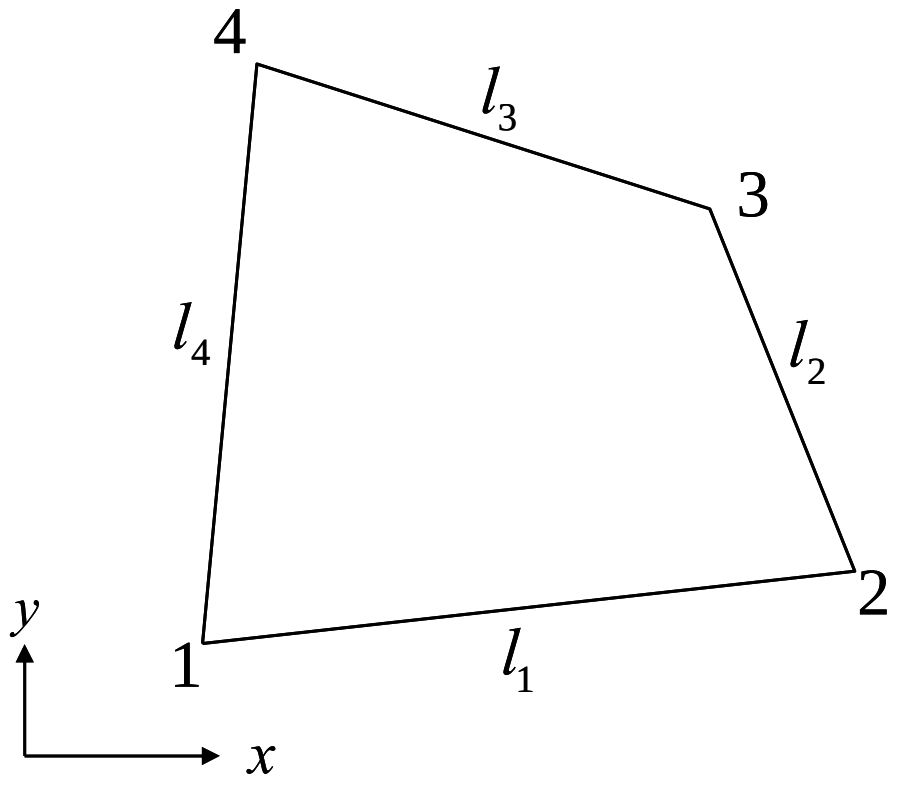

Figure 4.   A sample quadrilateral element

The wedge function $w_I(x, y) = 0$ corresponding to the *I*-th node is defined by the following product if the two edges do not pass through node *I*

$$w_I(x, y) = \kappa_I l_{I+1}(x, y) l_{I+2}(x, y) \tag{25}$$

where $\kappa_I$ are scalar constants. A compact definition of the Washspress shape function corresponding to the *I*-th node is given by

$$N_I(x, y) = \frac{w_I(x, y)}{\sum_J w_J(x, y)} \tag{26}$$

In this definition, a necessary requirement for the shape function is to be linear along the element boundary, which can be achieved by an appropriate selection of constants $\kappa_I$ as described in Ref [6].

It should be noted that the shape function in (26) is not always linear within the element, nor are the Lagrange shape functions used in the original smoothed FEM [21] because of the bilinear term '*xy*' in the basis. In the standard finite elements, an element can be divided into appropriate quadrilateral smoothing cells to avoid the nonlinearity of shape functions along the cell edges. In XFEM, however, the elements are generally cut by discontinuous interface with irregular geometries. This can be illustrated in Figure 5, in which a triangulation scheme used in standard XFEM is adopted to build the smoothing cells and the values of the shape function are plotted along the edges of smoothing cells. It is clearly seen that neither the standard nor the enriched shape functions are linear along the edges of smoothing cells in an enriched element. In these cases, integration of both the standard and the enrichment shape functions should be treated carefully to maintain accuracy and convergence rates. In this work, we will use three Gauss points in the enriched elements. In the standard elements, one Gauss point is enough to calculate the contour integral along each edge due to the linearity of the shape functions.

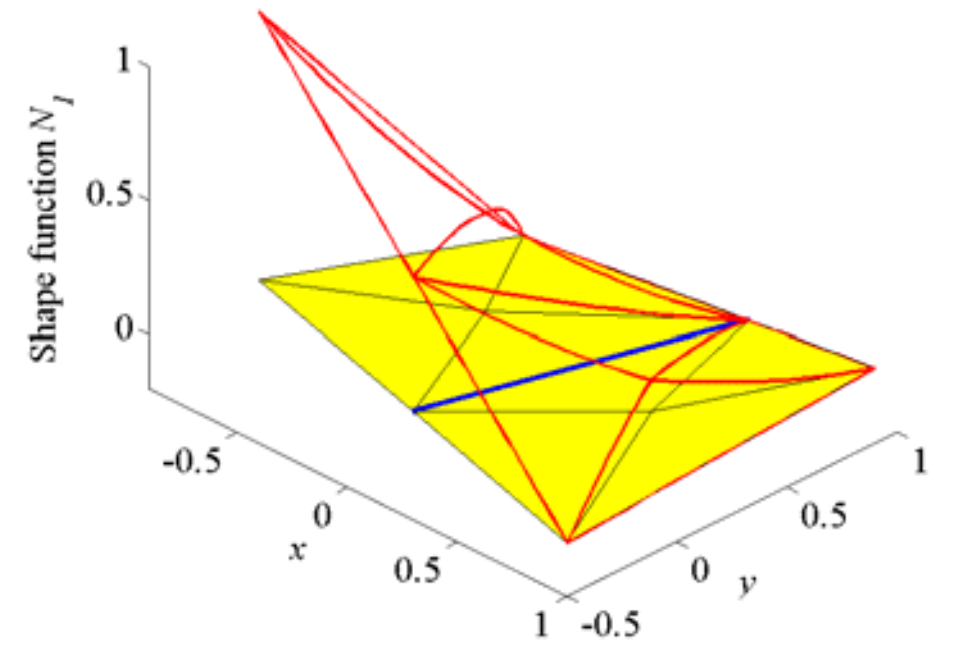

(a)

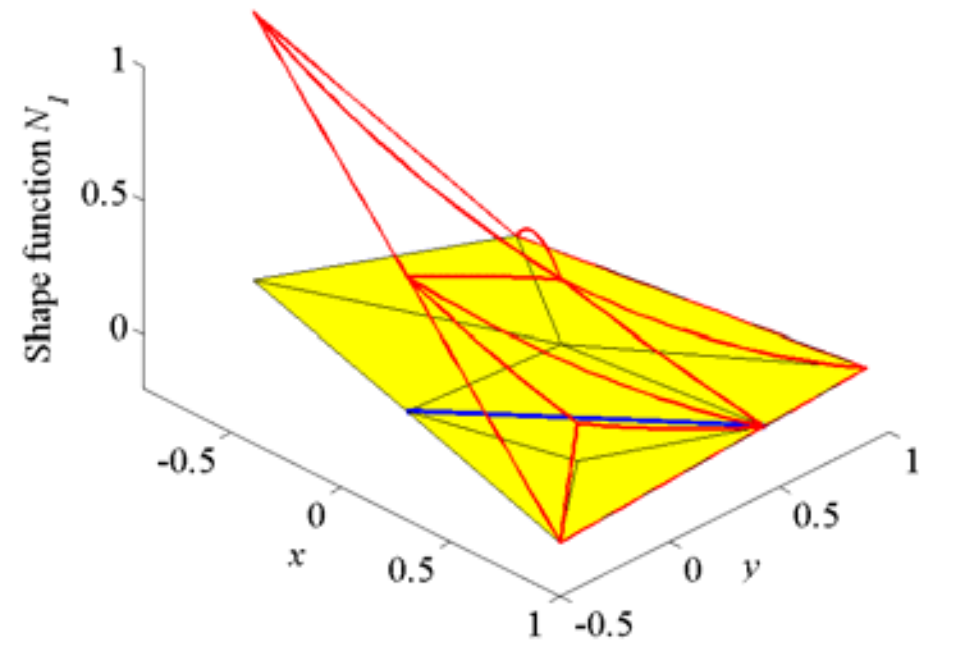

(b)

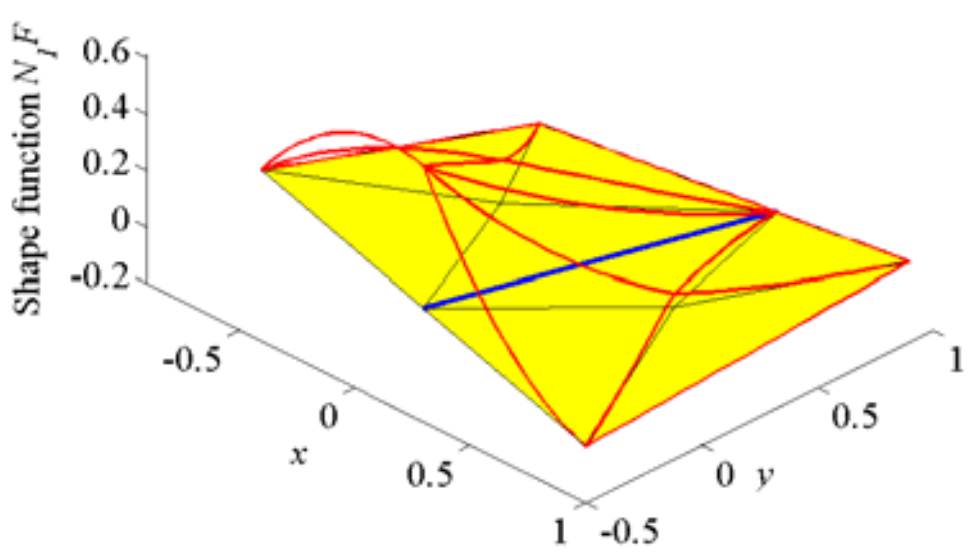

(c)

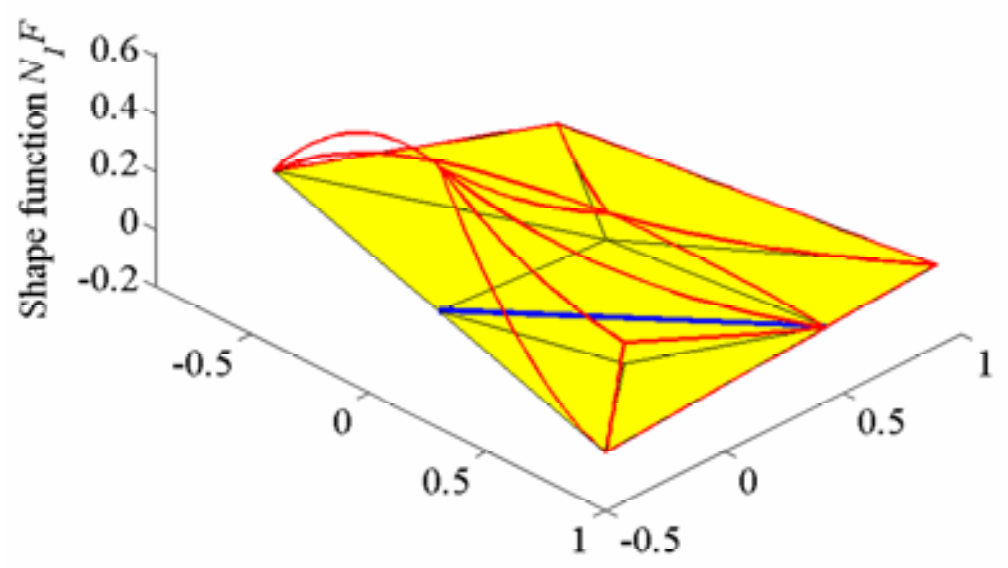

(d)

Figure 5. Triangular smoothing cell subdivision and shape functions along the edges of smoothing cells in an enriched element. (a) and (b): standard shape functions; (c) and (d): enriched shape functions. The blue bold lines represent the interfaces.

To evaluate the strain matrix, an element needs to be further partitioned into smoothing cells. It has been pointed out that when the number of smoothing cells in the element approaches infinity, the smoothed FEM solution will approach the standard FEM solution [22]. In our SmXFEM, we devised a new partition scheme of smoothing cells for the enriched elements.

For the quadrilateral elements, an element is firstly divided into $m \times n$ quadrilateral smoothing cells, where $m$ and $n$ are the number of subcells in the $x$- and $y$- directions, respectively. In the enriched elements, those smoothing cells are further distinguished as standard smoothing cells and enriched smoothing cells. Then the enriched smoothing cells are partitioned into a group of sub-triangles. This scheme is illustrated in Figure 6.

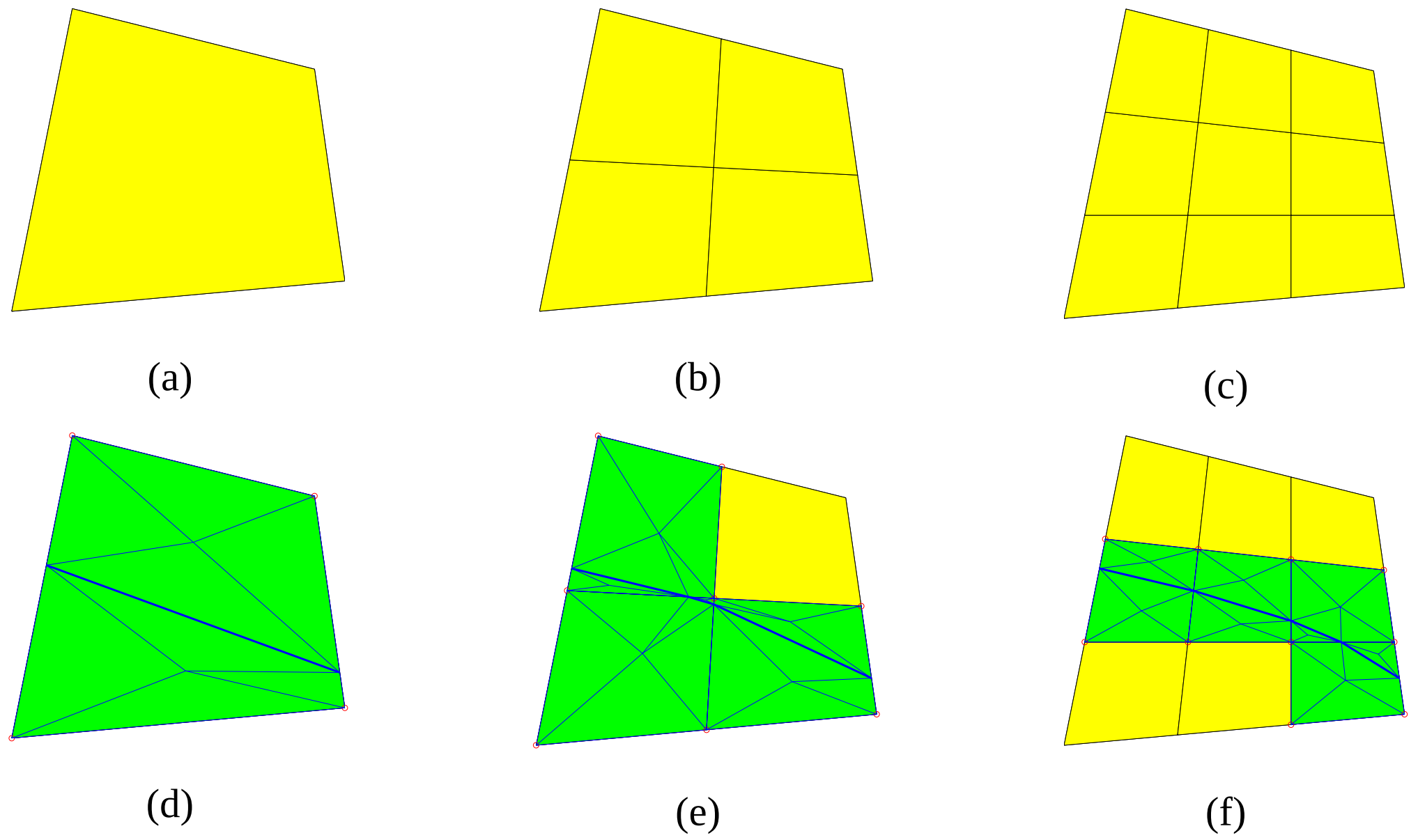

Figure 6. Partition of smoothing cells in an element. For the standard element: (a) $1 \times 1$ cell partition; (b) $2 \times 2$ cell partition; (c) $3 \times 3$ cell partition; For the enriched elements: (d) $1 \times 1$ cell partition; (e) $2 \times 2$ cell partition; (f) $3 \times 3$ cell partition. The yellow cells are standard smoothing cells and the green cells represent enriched smoothing cells.

It can be seen that the present partition scheme of smoothing cells provides great flexibility for arbitrary $m \times n$ cells partition according to the requirement of numerical accuracy and it is also extensible to higher order elements and three dimensional hexahedral elements. Another advantage of this scheme is that the curved interface can be more accurately resolved when the number of smoothing cells increases, which can be clearly seen in Figure 6(d)-(f).

## 5. Numerical examples

The SmXFEM developed above is suitable for a variety of problems dealing with nanoscale materials and structures with interface stress effects. In this paper, we focus on analyzing nano-inhomogeneities and their equilibrium shapes.

### *5.1 Convergence and computational efficiency*

Consider an isotropic circular inclusion with a coherent interface embedded in an isotropic elastic matrix of infinite extent, subjected to a dilatational eigenstrain $\boldsymbol{\varepsilon}^* = \varepsilon^* \mathbf{I}$ as shown in Figure 7. The analytical solution to this problem is given by [33]. The bulk elastic constants for aluminum are $\lambda$ = 58.17 GPa, $\mu$ = 26.13 GPa [23] and the interface elastic constants are set as $\tau^S = -1\,\text{N/m}$ and $L^S = 10\,\text{N/m}$ [40]. The radius of the inclusion is taken as $R$ = 5 nm. The eigenstrain $\varepsilon^*$ = 0.01 is used. The simulation is performed in a square domain of 15×15 nm and the boundary conditions are applied by imposing the exact displacements on the external boundary of the domain.

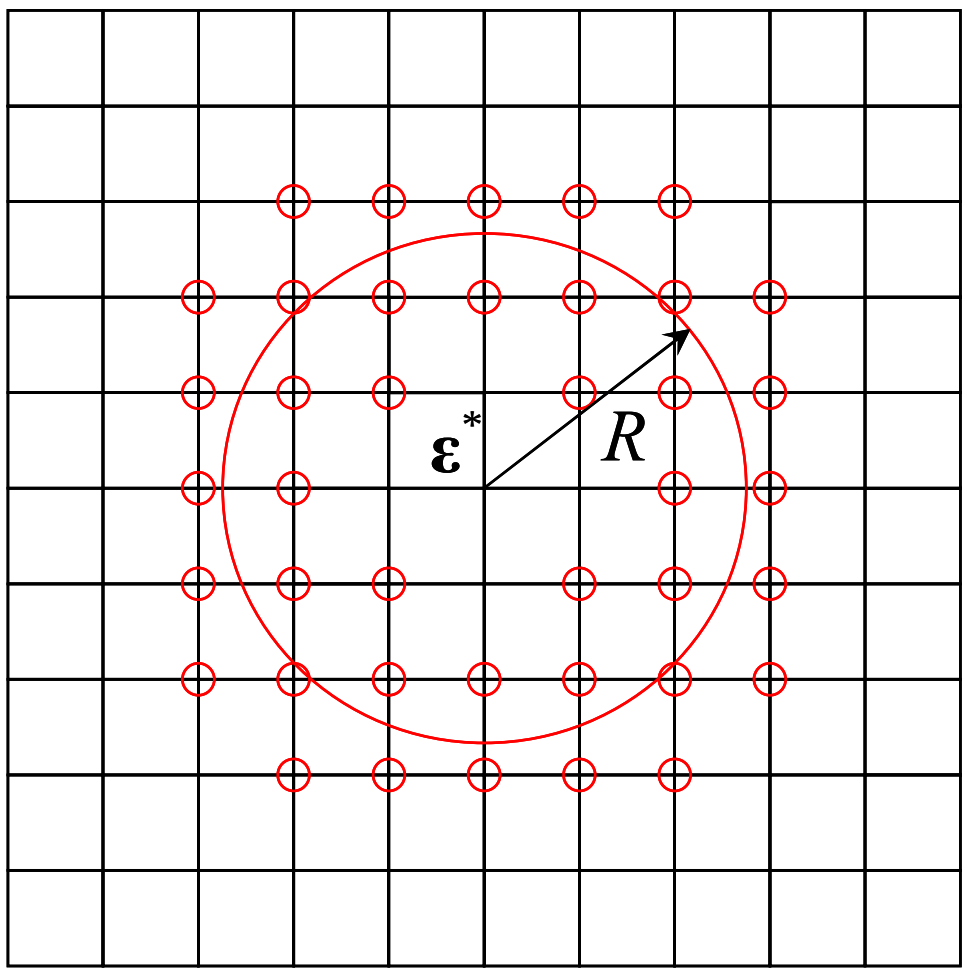

Figure 7. Computational model in smoothed XFEM with mesh, interface and enriched nodes.

In order to study the convergence of the present method, we define the following two norms, i.e. displacement norm $e_d$ and energy norm $e_e$,

$$e_d = \frac{\left\| \mathbf{u}^h - \mathbf{u}^{exact} \right\|}{\left\| \mathbf{u}^{exact} \right\|} \tag{27}$$

$$e_e = \sqrt{\frac{\int_\Omega \left( \boldsymbol{\varepsilon}^h(\mathbf{x}) - \boldsymbol{\varepsilon}^{exact}(\mathbf{x}) \right) : \mathbf{L} : \left( \boldsymbol{\varepsilon}^h(\mathbf{x}) - \boldsymbol{\varepsilon}^{exact}(\mathbf{x}) \right) d\Omega}{\int_\Omega \boldsymbol{\varepsilon}^{exact}(\mathbf{x}) : \mathbf{L} : \boldsymbol{\varepsilon}^{exact}(\mathbf{x}) d\Omega}} \tag{28}$$

where $\mathbf{u}^h$ and $\boldsymbol{\varepsilon}^h$ are the displacement and strain fields obtained from the present SmXFEM by using different element sizes $h$. The quantities with superscript '*exact*' are the exact solutions.

Figure 8 and Figure 9 show the displacement and energy norms for the cases without and with the interface stress, respectively. We can see that the optimal convergence rates are obtained for both smoothed and standard XFEM when the interfacial energy effects are ignored. The convergence rate is suboptimal when interfacial energy effects are considered and nearly 0.88 for the energy norm except for the 1×1 cell subdivision. The convergence rate is however much higher than that of 0.4-0.5 reported in [42] for the conventional XFEM. It is seen that the SmXFEM shows better performance over the standard XFEM if more than 1×1 cells are used. This is more obvious in Figure 9 where the interfacial energy effects are considered. We can also find that the relative errors decrease with the increase of the number of smoothing cells. However, when the number of smoothing cells exceeds 2×2, the improvement of the accuracy is not very significant.

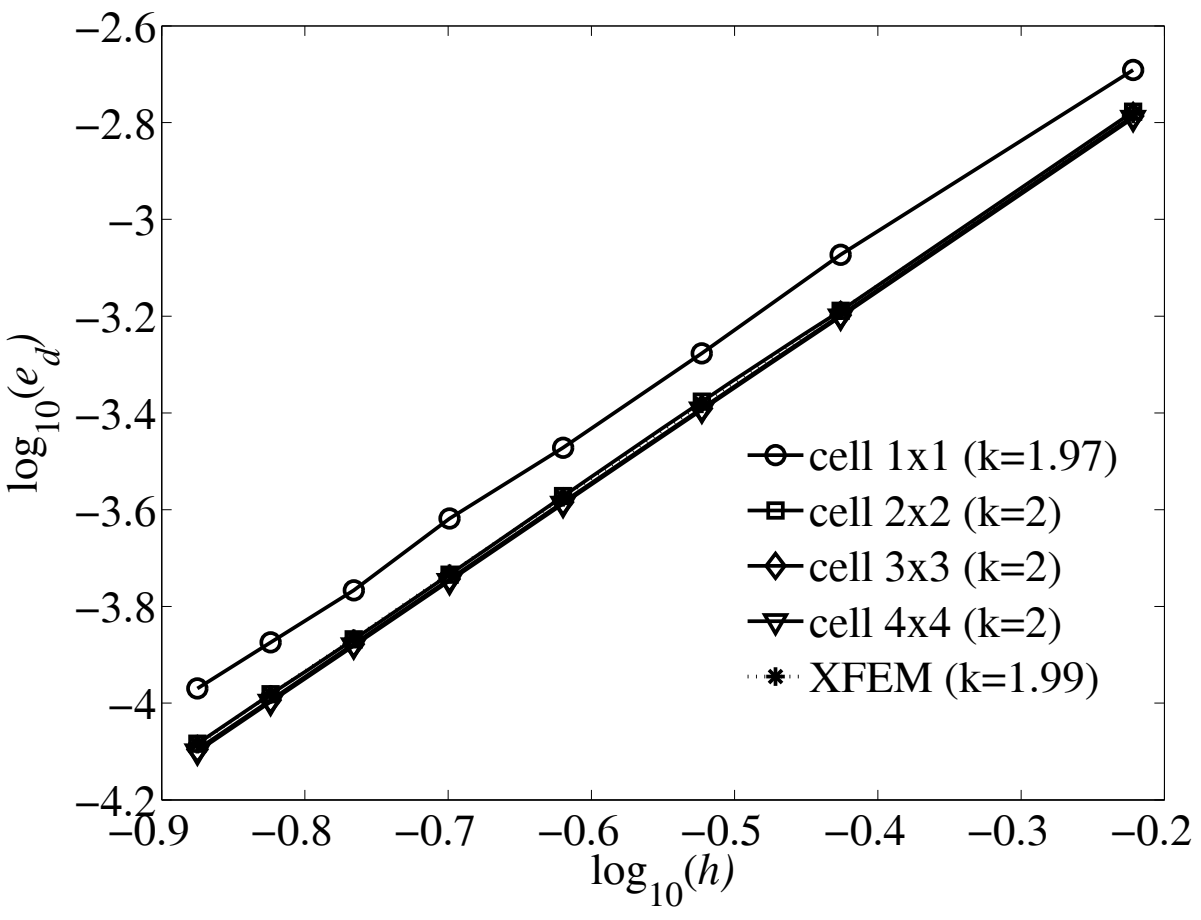

(a)

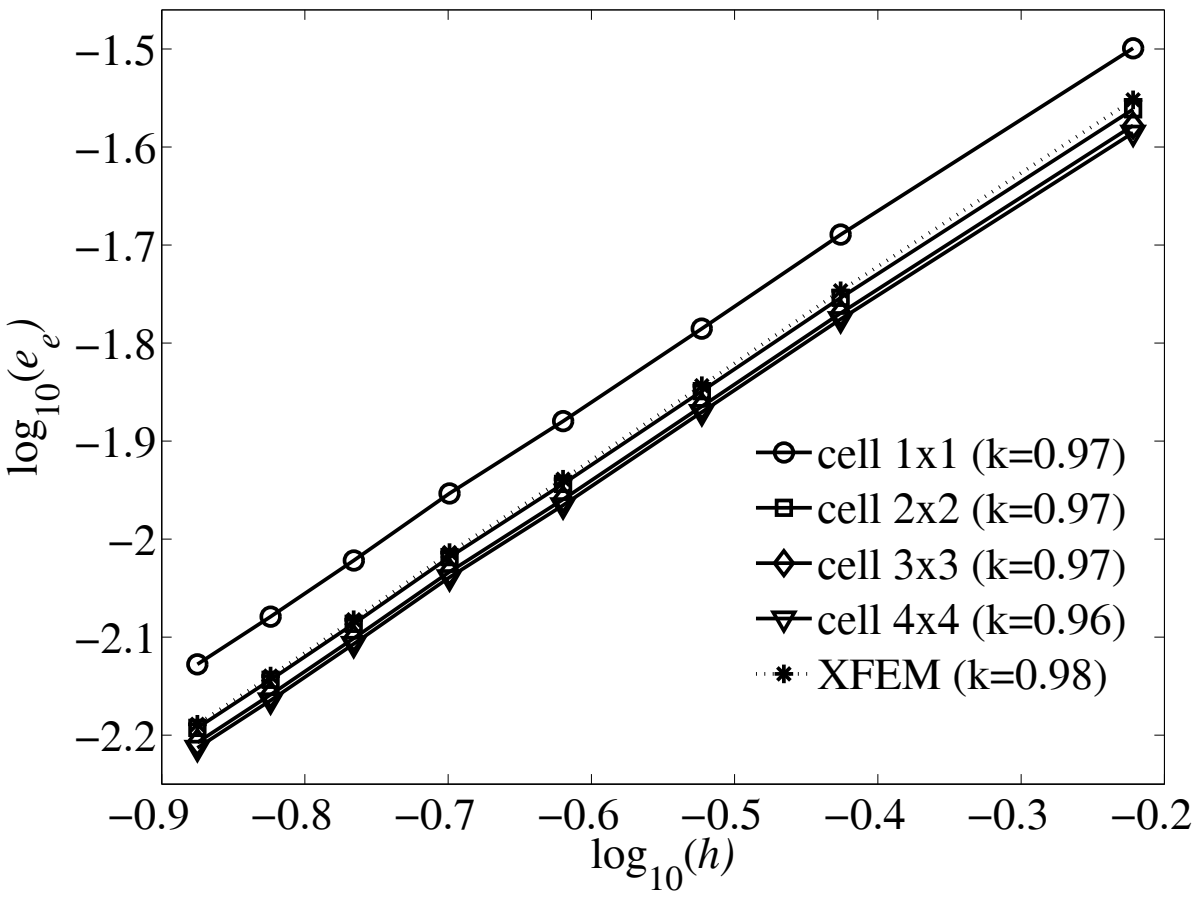

(b)

Figure 8. Convergence rate of circular inclusion problem without interface stress. (a) Displacement norm; (b) Energy norm.

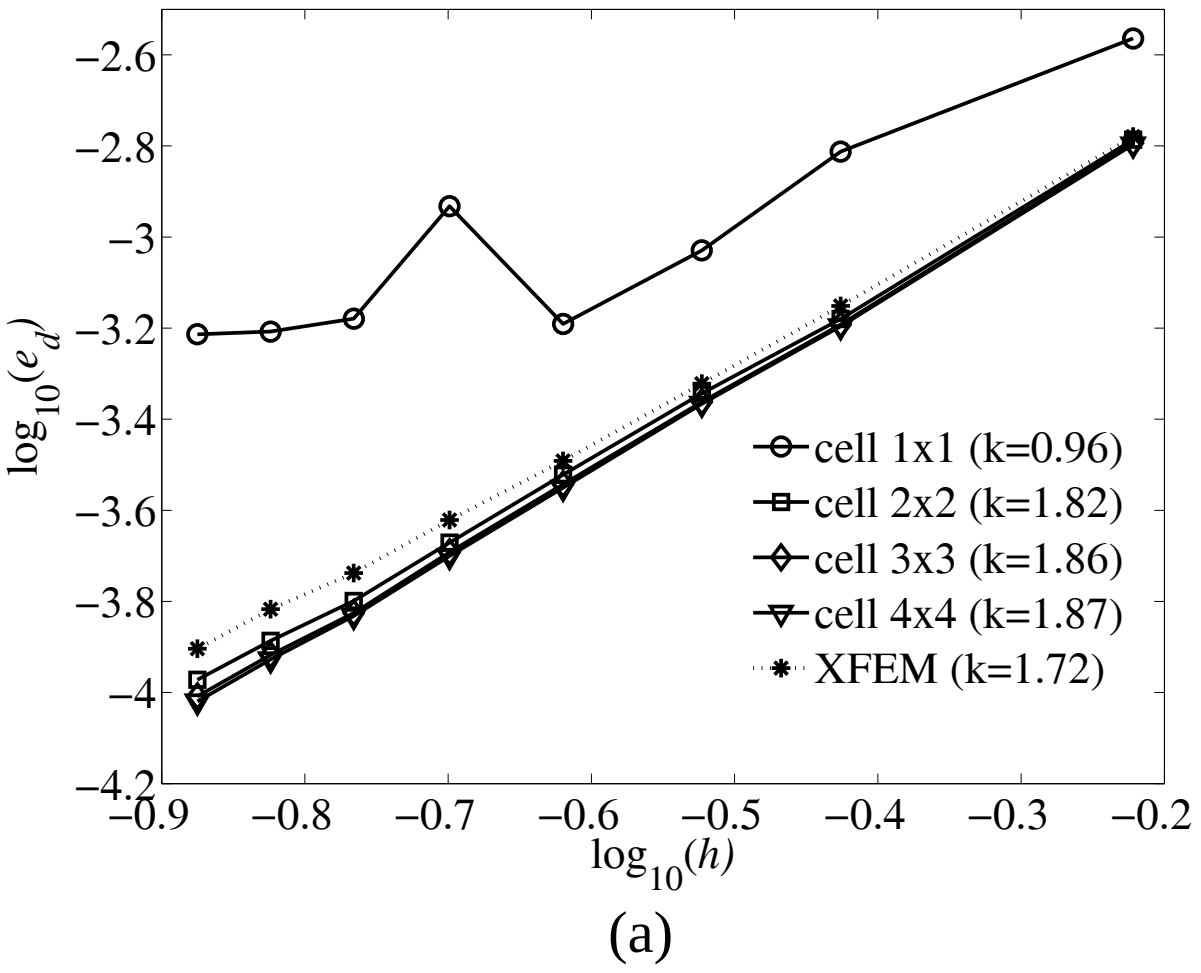

(a)

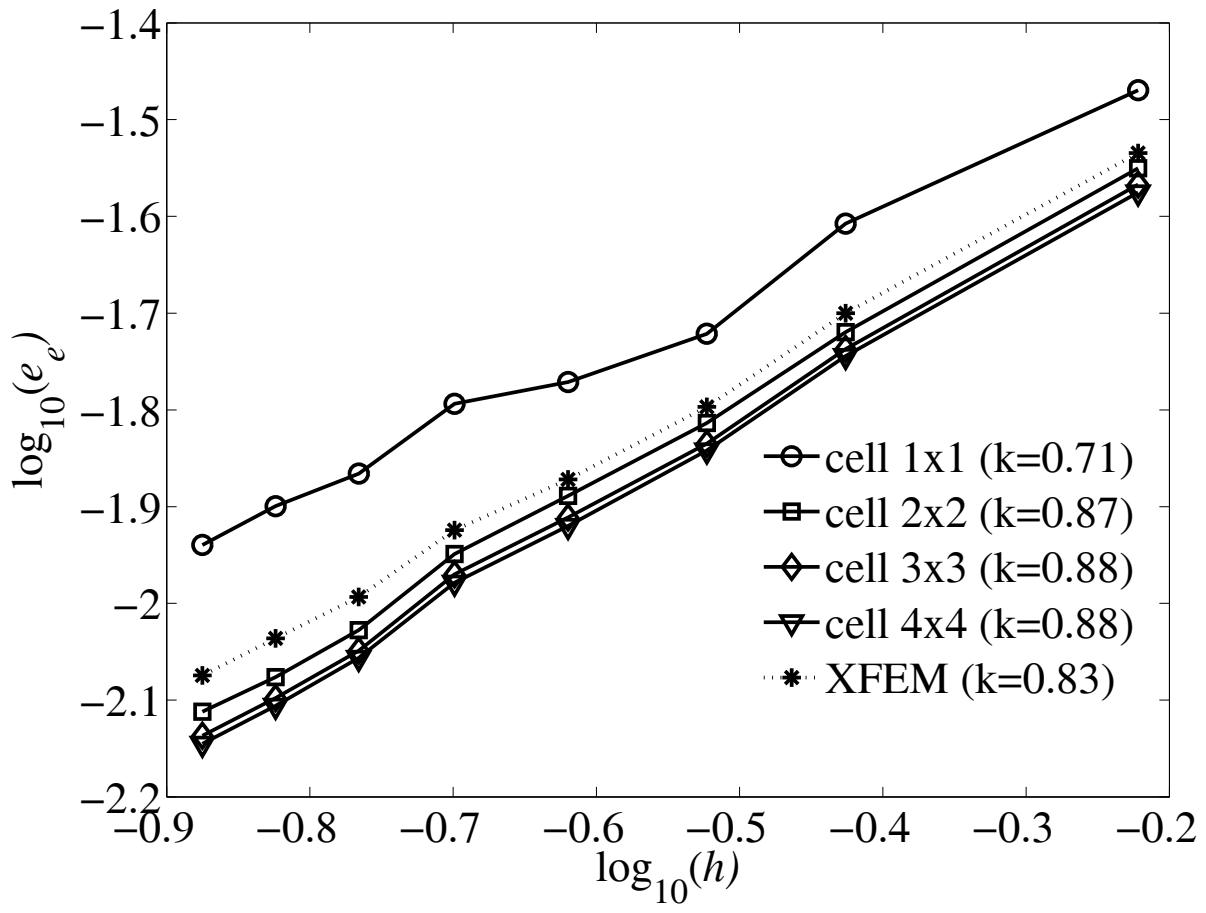

(b)

Figure 9. Convergence rate of circular inclusion problem with interface stress. (a) Displacement norm; (b) Energy norm.

As we mentioned, in the SmXFEM, the evaluation of the derivatives of shape functions and the Jacobian matrix in each element is not required to calculate the strain matrix. Moreover, the domain integrals are transformed into boundary integrals. In addition, in the standard XFEM, the numerical quadrature has to be performed over all the subtriangles in an enriched element. Therefore, two isoparametric mapping steps are required: mapping from the reference coordinate $\xi_{\Delta}$ of subtriangles to the global coordinate $\mathbf{x}$ and an inverse map from the global to reference coordinate system of the bi-unit square $\xi_{\square} = \mathbf{x}^{-1}\left(\xi_{\square}\right)$. Generally iterative algorithms are required for the second step mapping [38], which would be time consuming in the standard XFEM when the elements cut by the interface are numerous. However, the Wachspress shape function in our SmXFEM is defined directly on the global coordinate system without requiring any inverse coordinate mapping. Those characteristics can effectively reduce the computational cost.

Figure 10 compares the computational time versus total element numbers for the SmXFEM and the standard XFEM. The program runs on Windows XP operating system with Intel Xeon CPU 2.27 GHz. Since the only difference between XFEM and SmXFEM is the way to calculate the element stiffness matrix, we use the total time consumed to form the global matrix as the efficiency measure. It is clearly seen that the slope for SmXFEM is close to 1, which implies that the computational time is linearly proportional to the total number of elements. While for XFEM, this slope is nearly 1.86, which means the computational cost increases exponentially with the element numbers. This advantage makes SmXFEM very attractive for large scale computations without loss of accuracy. Furthermore, to balance the computational cost and accuracy, we can flexibly choose different combination of cell partition schemes, such as 2×2 cells in the standard elements and 4×4 cells in the enriched elements, which will be used in the following numerical examples.

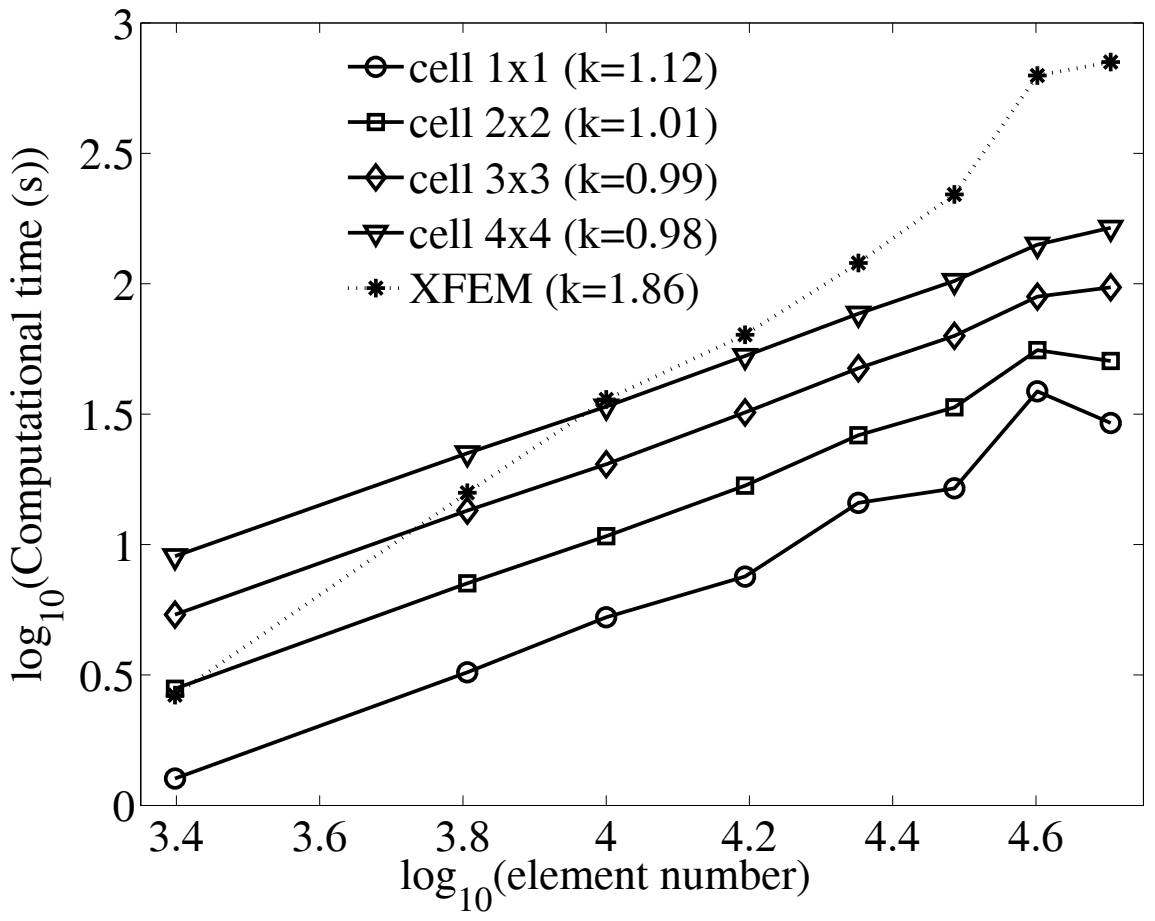

Figure 10. Computational time vs. element number.

5.2 *Elastic field of a nano-inclusion with interface effects*

The elastic field of the nano-inclusion problem described in section 5.1 is considered here. A regular mesh with 50×50 quadrilateral elements. The displacement and strain components in the radial direction with and without interface effects are compared and given in Figure 11. It is clearly seen that the present SmXFEM results are in very good agreement with the exact solutions, and capable of accurately capturing the interface effects. This again provides us with confidence on the validity of the approach to explore more complex problems.

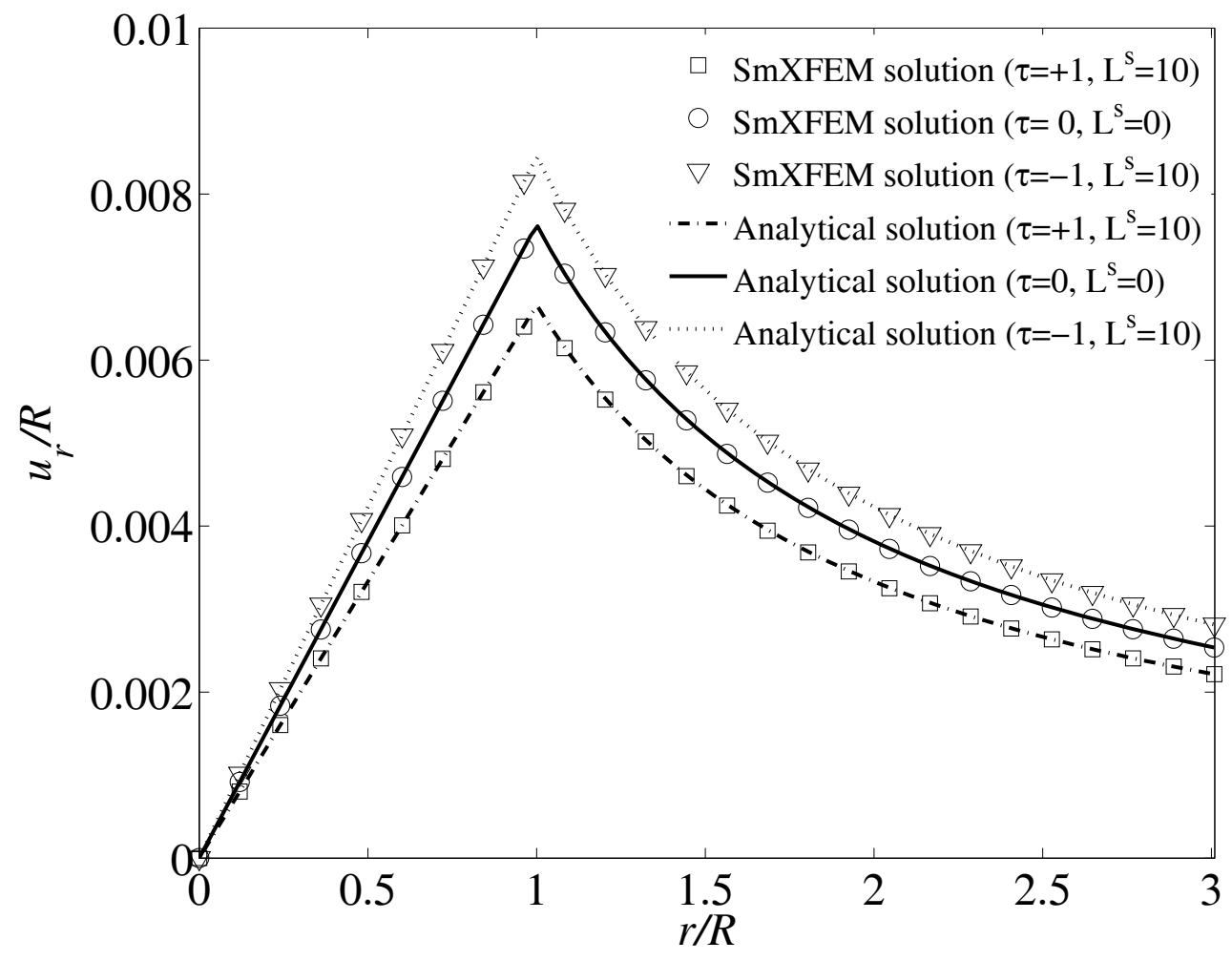

(a)

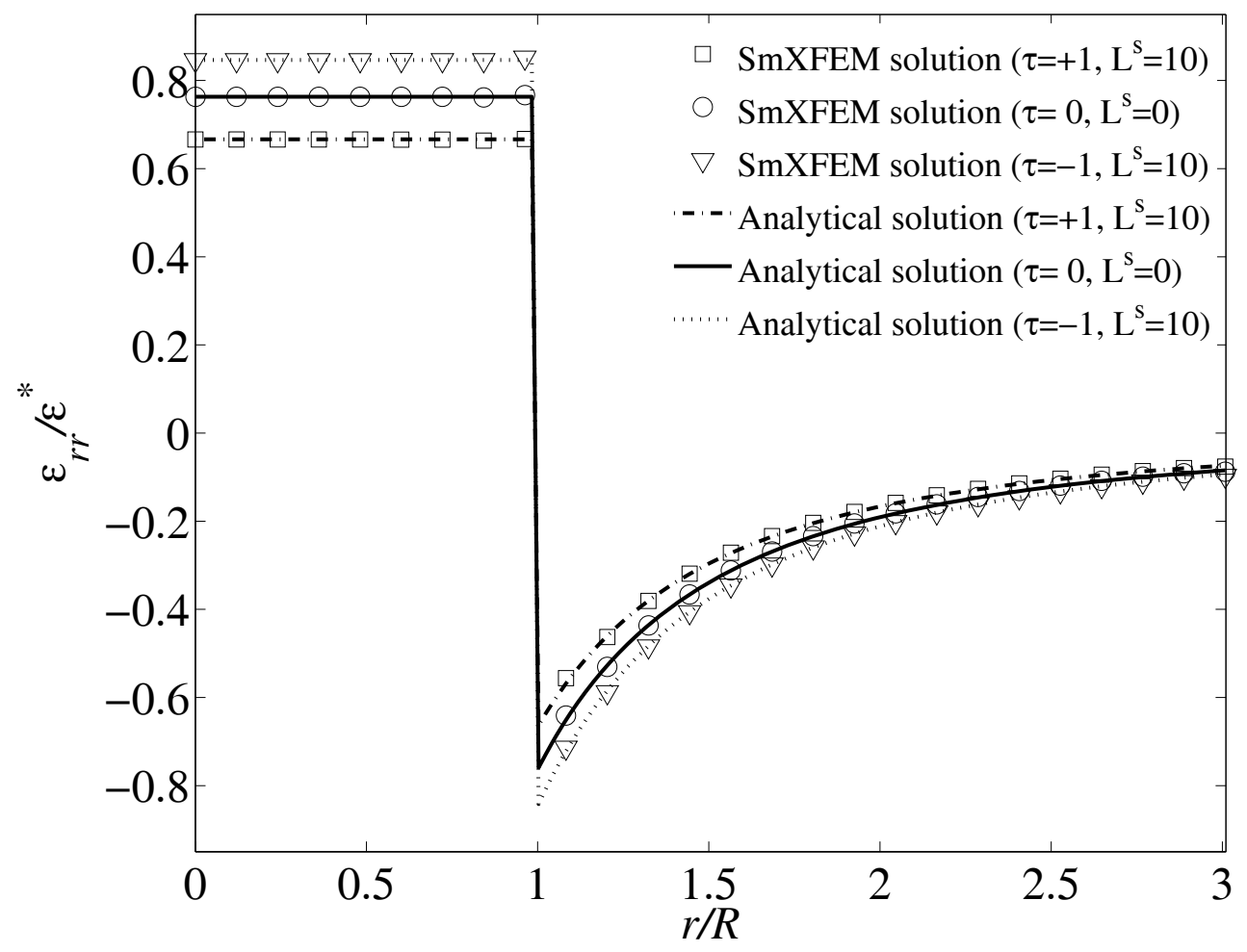

(b)

Figure 11. Comparison of the normalized radial displacement $u_r$ (a) and strain $\varepsilon_{rr}$ (b) from the exact and the present SmXFEM solutions.

*5.3 Equilibrium shapes of nano-particles in elastic solids*

It is experimentally observed that particles (precipitates) in nickel base superalloys may evolve, with increasing particle size, from its initial spherical shape to a cuboidal shape with flat edges and round corners. The equilibrium morphology of such particles is governed by the competition between elastic strain energy due to the misfit strain and the interfacial energy, so that the problem is equivalent to finding a group of discrete values of level set that minimizes the total system energy $\Pi$.

The driving force on the interface can be represented by the configurational force, which has the following form [31, 19]

$$v_n^{ext} = \mathbf{n} \cdot [\![\mathbf{\Sigma}]\!] \cdot \mathbf{n} - \gamma\kappa + \lambda \tag{29}$$

where $[\![\ ]\!]$ represents the quantity jump across the interface. $\mathbf{\Sigma} = w\mathbf{I} - \boldsymbol{\sigma} \cdot \nabla \mathbf{u}$ is the Eshelby energy momentum tensor where $w$ denotes the elastic energy density. $\kappa$ represents the interfacial mean curvature and

$$\lambda = \frac{\int_\Gamma \left( \mathbf{n} \cdot [\![\mathbf{\Sigma}]\!] \cdot \mathbf{n} + \gamma\kappa \right) dS}{\int_\Gamma dS} \tag{30}$$

In this problem, the velocity of the interface can be computed by the SmXFEM and the motion of material interfaces is tracked by the level set evolution equation (12). It is noted that the velocity field in the level set equation is defined in the spatial domain, but the velocity given in (29) is only defined on the interface. In the following example, a second order fast marching method [32] is used to extend the interface velocity to the entire computational domain, and the third order accurate Hamilton-Jacobi WENO finite difference spatial discretization and third order TVD Runge-Kutta (RK) time discretization schemes are used to solve the hyperbolic level set equation [29].

We use the anisotropic elastic constants of Ni for the matrix, which are $C_{11}^M$ =246.5 GPa, $C_{12}^M$ =147.3 GPa and $C_{44}^M$ =124.7 GPa in Voigt notation. We adopt the assumption that the elastic constants of particles are proportional to that of matrix, $C_{ij}^I = \alpha C_{ij}^M$ [31], so that we have hard particles when $\alpha > 1$, soft particles when $\alpha < 1$ and homogeneous particles when $\alpha = 1$. For simplicity, we further assume the interfacial energy is constant that $\gamma = \gamma_0$, where $\gamma_0$ =20 ergs cm$^{-2}$ and the dilatational misfit eigenstrain $\varepsilon^*$ = 0.3% [36]. Therefore, we can introduce the dimensionless characteristic length $L = C_{44}^{I}\varepsilon_0^2 R_0 / \gamma$, in which $R_0$ is the effective radius of particle, for example, $R_0 = \sqrt{A_0 / \pi}$ where $A_0$ is the area of a particle in two-dimensions. It is assumed that the initial shapes of the particles are circles for each case and the material axis <100> directions are the same as the coordinate directions.

Figure 12 shows the equilibrium shapes of an isolated particle at different sizes and stiffness ratios. It is seen that for the hard and homogeneous particles that $\alpha \geq 1$, the equilibrium morphologies undergo a transition from the circle-like to the square-like shapes with round corners when the particle size increases. This is in good agreement with experimental observations [11]. For the soft particles, the cuboidal shape remains for $L = 5$. However, the equilibrium shape gradually becomes concave with the increase of the particle size as shown in the uppermost row of Figure 12.

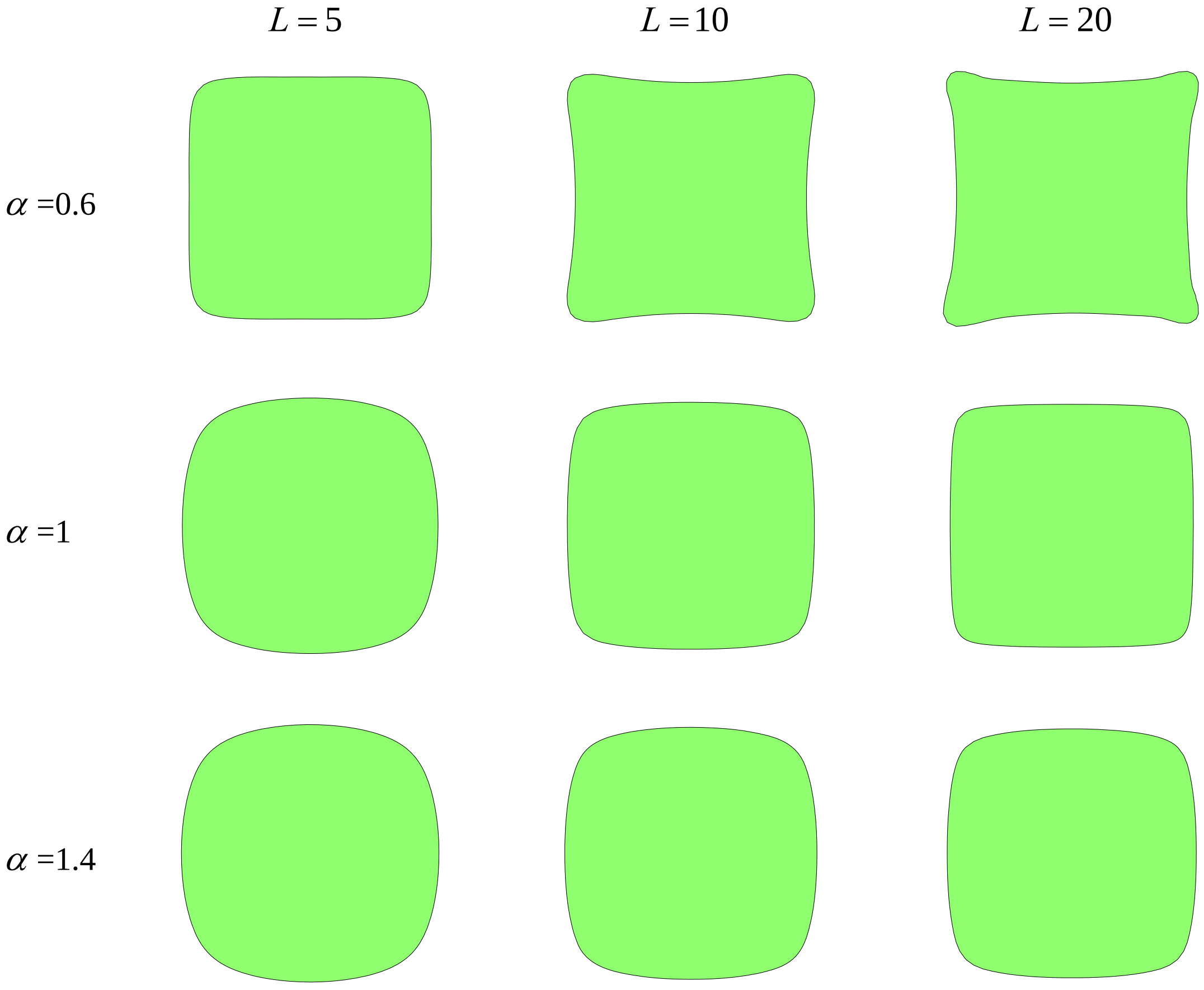

Figure 12. Equilibrium shapes of an isolated inhomogeneity for $\alpha$ =0.6, 1, 1.4 and $L$ =5, 10, 20.

It should be mentioned that the present results are very similar to that obtained from boundary integral methods [39, 31] and finite element methods [18], in which a number of marker points have to be placed on the geometric interface to track its motion. However, the proposed hybrid SmXFEM/LSM approach inherits the merits of both the extended finite element and level set method and provides distinct advantages to model the moving interface in solid materials. The main advantage of the SmXFEM is that the meshes do not have to conform to evolving interfaces, so that no remeshing is required during the evolution.

The implicit description of moving interfaces using the level set method naturally allows topological changes of phase boundaries, which are generally very hard to handle in the conventional marker particle method.

**6. Conclusion**

In this paper, we developed a hybrid smoothed extended finite element/level set method with evolving interfaces for modeling equilibrium shapes of nano-inhomogeneities. A linear interface elasticity model was adopted to account for interfacial excess energy effect. It is shown that once the interfacial energy is considered, the energetically favorable shapes of precipitate particles would depend on the particle size, misfit strain and the elastic constants of bulk materials due to the competition between the interfacial energy and anisotropic elastic energy. The SmXFEM developed here has several advantages over the conventional XFEM methods. For example, the Jacobian matrix and the derivatives of the standard and enrichment shape functions need not be computed and the domain integrals in the finite elements are replaced by boundary integrals, which can efficiently reduce the computational cost without loss of accuracy. The Wachspress shape functions are for the first time introduced into the SmXFEM framework, which greatly facilitates accurate numerical integration along the edges of smoothing cells. The accuracy and convergence rates are studied in the numerical results, which show superior to the standard XFEM [42]. Further, we showed that the computational time of our SmXFEM scales linearly with the number of elements, while that of our implementation of the XFEM increases exponentially (exponent 1.86) with the number of elements. We believe that the hybrid SmXFEM/LSM developed here is an effective numerical tool for nano-scale material analysis and designs, and can be used to study nano-scale inhomogeneous materials and structures with complex morphological changes.

**Reference:**

[1] Bordas S, Nguyen P V, Dunant C, Guidoum A and Nguyen-Dang H 2007 An extended finite element library *Int J Numer Meth Eng* **71** 703-32

[2] Bordas S P A and Natarajan S 2010 On the approximation in the smoothed finite element method (SFEM) *Int J Numer Meth Eng* **81** 660-70

[3] Chen J S, Wu C T, Yoon S and You Y 2001 A stabilized conforming nodal integration for Galerkin mesh-free methods *Int J Numer Meth Eng* **50** 435-66

[4] Chessa J and Belytschko T 2003 An extended finite element method for two-phase fluids *J Appl Mech-T Asme* **70** 10-7

[5] Chessa J, Smolinski P and Belytschko T 2002 The extended finite element method (XFEM) for solidification problems *Int J Numer Meth Eng* **53** 1959-77

[6] Dasgupta G 2003 Interpolants within convex polygons: Wachspress' shape functions *J Aerospace Eng* **16** 1-8

[7] Dingreville R, Qu J M and Cherkaoui M 2005 Surface free energy and its effect on the elastic behavior of nano-sized particles, wires and films *J Mech Phys Solids* **53** 1827-54

[8] Duan H L, Wang J, Huang Z P and Karihaloo B L 2005 Size-dependent effective elastic constants of solids containing nano-inhomogeneities with interface stress *J Mech Phys Solids* **53** 1574-96

[9] Duddu R, Bordas S, Chopp D and Moran B 2008 A combined extended finite element and level set method for biofilm growth *Int J Numer Meth Eng* **74** 848-70

[10] Duddu R, Chopp D L and Moran B 2009 A Two-Dimensional Continuum Model of Biofilm Growth Incorporating Fluid Flow and Shear Stress Based Detachment *Biotechnol Bioeng* **103** 92-104

[11] Fahrmann M, Fratzl P, Paris O, Fahrmann E and Johnson W C 1995 Influence of Coherency Stress on Microstructural Evolution in Model Ni-Al-Mo Alloys *Acta Metall Mater* **43** 1007-22

[12] Farsad M, Vernerey F J and Park H S 2010 An extended finite element/level set method to study surface effects on the mechanical behavior and properties of nanomaterials *Int J Numer Meth Eng* **84** 1466-89

[13] Fries T P and Belytschko T 2010 The extended/generalized finite element method: An overview of the method and its applications *Int J Numer Meth Eng* **84** 253-304

[14] Gao W, Yu S W and Huang G Y 2006 Finite element characterization of the size-dependent mechanical behaviour in nanosystems *Nanotechnology* **17** 1118-22

[15] Gracie R, Oswald J and Belytschko T 2008 On a new extended finite element method for dislocations: Core enrichment and nonlinear formulation *J Mech Phys Solids* **56** 200-14

[16] Gurtin M E and Murdoch A I 1975 Continuum Theory of Elastic-Material Surfaces *Archive for Rational Mechanics and Analysis* **57** 291-323

[17] Gurtin M E, Weissmuller J and Larche F 1998 A general theory of curved deformable interfaces in solids at equilibrium *Philos Mag A* **78** 1093-109

[18] Jog C S, Sankarasubramanian R and Abinandanan T A 2000 Symmetry-breaking transitions in equilibrium shapes of coherent precipitates *J Mech Phys Solids* **48** 2363-89

[19] Kolling S, Mueller R and Gross D 2003 The influence of elastic constants on the shape of an inclusion *Int J Solids Struct* **40** 4399-416

[20] Legay A, Chessa J and Belytschko T 2006 An Eulerian-Lagrangian method for fluid-structure interaction based on level sets *Comput Method Appl M* **195** 2070-87

[21] Liu G R, Dai K Y and Nguyen T T 2007 A smoothed finite element method for mechanics

problems *Comput Mech* **39** 859-77

[22] Liu G R, Nguyen T T, Dai K Y and Lam K Y 2007 Theoretical aspects of the smoothed finite element method (SFEM) *Int J Numer Meth Eng* **71** 902-30

[23] Meyers M A and Chawla K K 2009 *Mechanical behavior of materials* (Cambridge, UK ; New York: Cambridge University Press)

[24] Miller R E and Shenoy V B 2000 Size-dependent elastic properties of nanosized structural elements *Nanotechnology* **11** 139-47

[25] Moës N, Cloirec M, Cartraud P and Remacle J F 2003 A computational approach to handle complex microstructure geometries *Comput Method Appl M* **192** 3163-77

[26] Moës N, Dolbow J and Belytschko T 1999 A finite element method for crack growth without remeshing *Int J Numer Meth Eng* **46** 131-50

[27] Moumnassi M, Belouettar S, Bechet E, Bordas S P A, Quoirin D and Potier-Ferry M 2011 Finite element analysis on implicitly defined domains: An accurate representation based on arbitrary parametric surfaces *Comput Method Appl M* **200** 774-96

[28] Olsson P A T and Park H S 2012 On the importance of surface elastic contributions to the flexural rigidity of nanowires *J Mech Phys Solids* **60** 2064-83

[29] Osher S and Fedkiw R P 2003 *Level set methods and dynamic implicit surfaces* (New York: Springer)

[30] Osher S and Sethian J A 1988 Fronts Propagating with Curvature-Dependent Speed - Algorithms Based on Hamilton-Jacobi Formulations *J Comput Phys* **79** 12-49

[31] Schmidt I and Gross D 1997 The equilibrium shape of an elastically inhomogeneous inclusion *J Mech Phys Solids* **45** 1521-49

[32] Sethian J A 1999 Fast marching methods *Siam Rev* **41** 199-235

[33] Sharma P and Ganti S 2004 Size-dependent Eshelby's tensor for embedded nano-inclusions incorporating surface/interface energies *J Appl Mech-T Asme* **71** 663-71

[34] She H and Wang B A 2009 A geometrically nonlinear finite element model of nanomaterials with consideration of surface effects *Finite Elements in Analysis and Design* **45** 463-7

[35] Shuttleworth R 1950 The Surface Tension of Solids *P Phys Soc Lond A* **63** 444-57

[36] Su C H and Voorhees P W 1996 The dynamics of precipitate evolution in elastically stressed solids .1. Inverse coarsening *Acta Mater* **44** 1987-99

[37] Sukumar N, Chopp D L, Moes N and Belytschko T 2001 Modeling holes and inclusions by level sets in the extended finite-element method *Comput Method Appl M* **190** 6183-200

[38] Sukumar N and Prevost J H 2003 Modeling quasi-static crack growth with the extended finite element method Part I: Computer implementation *Int J Solids Struct* **40** 7513-37

[39] Thompson M E, Su C S and Voorhees P W 1994 The Equilibrium Shape of a Misfitting Precipitate *Acta Metall Mater* **42** 2107-22

[40] Tian L and Rajapakse R K N D 2007 Finite element modelling of nanoscale inhomogeneities in an elastic matrix *Comp Mater Sci* **41** 44-53

[41] Wachspress E L 1975 *A rational finite element basis* (New York: Academic Press)

[42] Yvonnet J, Le Quang H and He Q C 2008 An XFEM/level set approach to modelling surface/interface effects and to computing the size-dependent effective properties of nanocomposites *Comput Mech* **42** 119-31

[43] Zhao X J and Qu J 2012 Effects of interfacial excess energy on the elastic field of a nano-inhomogeneity *Mech Mater* **55** 41-8